\input vanilla.sty
\scaletype{\magstep1}
\scalelinespacing{\magstep1}
\def\bull{\vrule height .9ex width .8ex depth -.1ex}

\title Surjective isometries on rearrangement-invariant
spaces \endtitle

\author N.J.  Kalton and Beata
Randrianantoanina\footnote{Both authors were supported by
NSF-grant DMS-9201357\newline This research will form part
of the Ph.D. thesis of the second author currently under
preparation at the University of Misssouri\newline AMS
Classification 46B04}\\ Department of Mathematics\\
University of Missouri\\ Columbia\\ Mo. 65211 \endauthor

\vskip2truecm

\subheading{Abstract}We prove that if $X$ is a real
rearrangement-invariant function space on $[0,1]$, which is
not isometrically isomorphic to $L_2,$ then every surjective
isometry $T:X\to X$ is of the form $Tf(s)=a(s)f(\sigma(s))$
for a Borel function $a$ and an invertible Borel map
$\sigma:[0,1] \to [0,1].$ If $X$ is not equal to $L_p$, up
to renorming, for some $1\le p\le \infty$ then in addition
$|a|=1$ a.e. and $\sigma$ must be measure-preserving.

\vskip2truecm

\subheading{1.  Introduction}

The main result of this paper is the following theorem,
which combines the statements of Theorems 6.4 and 7.2 below.
We denote Lebesgue measure on $[0,1]$ by $\lambda$ and use
the term rearrangement-invariant Banach function space in
the sense of Lindenstrauss-Tzafriri [19].

\proclaim{Theorem 1.1}Let $X$ be a (real)
rearrangement-invariant Banach function space on $[0,1].$
Suppose $X$ is not (isometrically) equal to $L_2[0,1].$ Let
$T:X\to X$ be a surjective isometry.  Then \newline (1)
There exists a nonvanishing Borel function $a:[0,1]\to\bold
R$ and an invertible Borel map $\sigma:[0,1]\to[0,1]$ such
that, for any Borel set $B\subset [0,1],$ we have
$\lambda(\sigma^{-1}B)=0$ if and only if $\lambda(B)=0$ and
so that $Tf(s)=a(s)f(\sigma(s))$ a.e. for any $f\in X.$
\newline (2) If $X$ is not equal to $L_p$ for some $1\le
p\le \infty$ up to renorming, then $|a|=1$ a.e. and $\sigma$
is measure-preserving.  \endproclaim

The first part of this theorem is known for the case of
spaces of complex functions and is due to Zaidenberg [32],
[33]; we discuss the relationship between our result and
Zaidenberg's below.  The second part also holds for the
complex case, and is apparently new even in this case.

The study of isometries on classical function spaces goes
back to Banach [1] who proved that the isometries of
$L_p[0,1]$ are disjointness-preserving when $p\neq 2$ (see
[1] p. 175).  Lamperti [18] later characterized the
isometries on $L_p.$ Since then there has developed an
extensive literature on isometries of particular function
spaces; see [3], [4], [5] and [13] for example.

In the case of (not necessarily rearrangement-invariant)
complex function spaces, a technique developed by Lumer
[20], [21], [22] has proved particularly effective.  This
technique was used by Lumer [21], [22] to study isometries
on reflexive Orlicz spaces and later by Zaidenberg [32] and
[33] to study isometries on general r.i. spaces, $X$.  The
idea is to characterize first the hermitian operators
$H:X\to X.$ $H$ is hermitian if $\exp(itH)$ is an isometry
for every real $t.$ One shows that hermitian operators are
simply multiplication operators by real functions, unless
$X=L_2.$ Then, if $U$ is a surjective isometry on $X$ we
have that $UHU^{-1}$ is hermitian for every hermitian $H.$
Combining these ideas leads to Theorem 1.1 (1) in the
complex case.  See also, for example, [6], [7], [10] and
[31].  For a fuller discussion of the existing literature we
refer the reader to the forthcoming survey of Fleming and
Jamison [8].

This line of argument simply does not work for real spaces,
and most of the known results use geometric techniques (e.g.
extreme point arguments) for special spaces.  In this paper,
we follow a line of reasoning which is distantly related to
the Lumer technique.  We use the notion of a numerically
positive operator [26]; this is an operator $T$ such that
$\|\exp(-tT)\| \le 1$ for all $t\ge 0$ (see also [23] where
$-T$ is called dissipative).  Unfortunately, this is far too
weak a notion to allow us to characterize such operators on
an r.i. space, but by studying rank-one numerically positive
operators and using results of Flinn (see [26]) we are able
to prove a representation theorem for surjective isometries
(Proposition 6.3), which is a partial step towards our main
result.  Then by a probabilistic technique we obtain Theorem
6.4 which is equivalent to Theorem 1.1 (1).  Finally in
Theorem 7.2 we show that if $X$ is not $L_p$ up to renorming
then the representation in Theorem 6.4 can be further
narrowed to the trivial case as in Theorem 1.1 (2).

Some remarks on the nature of our results are in order.
First notice that we must restrict ourselves (as do Lumer
and Zaidenberg) to surjective isometries.  Many of the
special results quoted above apply equally to isometries
which are not surjective.  Secondly, it will be seen that
there appear to be obstacles to extending the main result to
r.i. spaces on $[0,\infty),$ (see, for example, [11], [12]).
However our results do apply equally to separable and
nonseparable r.i. spaces on $[0,1].$ The proof can be
simplified a little in the separable case and we indicate
such simplifications as various points in the paper.

\demo{Acknowledgement}We wish to thank Jim Jamison and Anna
Kaminska for bringing this question to our attention and for
providing copies of both [8] and a translation of
[33].\enddemo

\vskip2truecm

\subheading{2.  Introductory remarks on K\"othe function
spaces}

Let us suppose that $\Omega$ is a Polish space and that
$\mu$ is a $\sigma-$finite Borel measure on $\Omega.$ We use
the term K\"othe space in the sense of [19] p. 28.  Thus a
K\"othe function space $X$ on $(\Omega,\mu)$ is a Banach
space of (equivalence classes of) locally integrable Borel
functions $f$ on $\Omega$ such that:  \newline (1) If
$|f|\le |g|$ a.e. and $g\in X$ then $f\in X$ with $\|f\|_X
\le \|g\|_X.$\newline (2) If $A$ is a Borel set of finite
measure then $\chi_A\in X.$

We say that $X$ is order-continuous if whenever $f_n\in X$
with $f_n\downarrow 0$ a.e. then $\|f_n\|_X\downarrow 0.$
$X$ has the Fatou property if whenever $0\le f_n\in X$ with
$\sup\|f_n\|_X<\infty$ and $f_n\uparrow f$ a.e. then $f\in
X$ with $\|f\|_X=\sup\|f_n\|_X.$

The K\"othe dual of $X$ is denoted $X';$ Thus $X'$ is the
K\"othe space of all $g$ such that $\int
|f||g|\,d\mu<\infty$ for every $f\in X$ equipped with the
norm $\|g\|_{X'}= \sup_{\|f\|_X\le 1}\int |f||g|\,d\mu.$
Then $X'$ can be regarded as a closed subspace of the dual
$X^*$ of $X$.  If $X$ is order-continuous then $X'=X^*$; if
$X$ has the Fatou property then $X'$ is a norming subspace
of $X^*.$

A rearrangement-invariant function space (r.i. space) is a
K\"othe function space on $([0,1],\lambda)$ where $\lambda$
is Lebesgue measure which satisfies the conditions:\newline
(1) Either $X$ is order-continuous or $X$ has the the Fatou
property.  \newline (2) If $\tau:[0,1]\to [0,1]$ is any
measure-preserving invertible Borel automorphism then $f\in
X$ if and only if $f\circ\tau\in X$ and
$\|f\|_X=\|f\circ\tau\|_X.$ \newline (3)
$\|\chi_{[0,1]}\|_X=1.$

In this section we will make some introductory remarks about
operators and isometries on K\"othe function spaces.  Let us
suppose that $X$ is a K\"othe function space on
$(\Omega,\mu).$ We first consider those operators $T:X\to X$
which are continuous for the topology $\sigma(X,X').$ This
is equivalent to requiring the existence of adjoint
$T':X'\to X'.$ Of course, if $X$ is order-continuous, every
operator $T:X\to X$ is $\sigma(X,X')-$continuous.

Our first result is well-known, but we know no explicit
reference.

\proclaim{Proposition 2.1}Let $X$ be a K\"othe function
space on $(\Omega,\mu).$ The following conditions on $T:X\to
X$ are equivalent:  \newline (1) $T$is
$\sigma(X,X')-$continuous.\newline (2) If $0\le f_n\in X$
and $f_n\uparrow f$ a.e. then $\lim_{n\to\infty}\int
|h||Tf_n-Tf|d\mu=0$ for every $h\in X'.$ (3) If $0\le g\in
X$ and $|f_n|\le g$ with $f_n\to f$ a.e. then
$\lim_{n\to\infty}\int |h||Tf_n-Tf|d\mu=0$ for every $h\in
X'.$ \endproclaim

\demo{Remark}Note that (2) says that $T:X\to L_1(|h|d\mu)$
is an order-continuous operator for every $h\in X'$; see
Weis [29].\enddemo

\demo{Proof}$(1)\rightarrow (3):$ Consider the operator
$S:L_{\infty}(\mu)\to L_1(\mu)$ defined by $S\phi= h(T\phi
g).$ Then $S$ is $\sigma(L_{\infty},L_1)\to
\sigma(L_1,L_{\infty})$-continuous and hence weakly compact.
Now we may choose $\phi_n\to \phi$ a.e. such that
$\phi_ng=f_n$ and $\phi g=f.$ Consider the adjoint
$S':L_{\infty}\to L_1.$ Then $S'(B_{L_{\infty}})$ is weakly
compact and hence uniformly integrable in $L_1$ ([30], p.
137).  Thus $$\lim_{n\to\infty}\sup_{\|h\|_{\infty}\le
1}\left|\int \phi_n S'h\, d\mu\right|=0$$ which quickly
gives (3).

$(3)\rightarrow (2):$ Obvious.

$(2)\rightarrow (1):$ We must show that the adjoint
$T^*:X^*\to X^*$ maps $X'$ into $X'.$ Consider any $h\in
X'.$ It follows quickly that the set-function $A\to \int
hT\chi_Ad\mu$ is countably additive when restricted to any
Borel set $A$ of finite measure.  Thus there exists $\phi\in
L_0$ so that $\int h(Tf)d\mu=\int f\phi \,d\mu$ for any
simple function $f$ supported on a set of finite measure.
Now for general positive $f\in X$ we find a sequence $f_n$
of such simple functions so that $f_n\uparrow f$ a.e. and
then $$\int h(Tf)=\lim_{n\to\infty}\int
h(Tf_n)d\mu=\lim_{n\to\infty}\int \phi f_nd\mu=\int\phi
fd\mu.$$ We conclude that $T^*h=\phi\in X'.$ \bull\enddemo

An operator $T:X\to X$ will be called {\it elementary} if
there is a Borel function $a$ and a Borel map
$\sigma:\Omega\to\Omega$ such that $Tf(s)=a(s)f(\sigma(s))$
a.e. for every $f\in X.$ Observe that a necessary condition
on $a$ and $\sigma$ is that if $B$ is a Borel set with
$\mu(B)=0$ then $\mu(\sigma^{-1}B\cap\{|a|>0\})=0.$ $T$ is
called {\it disjointness-preserving} if $\min(|f|,|g|)=0$
a.e. implies $\min(|Tf|,|Tg|)=0,$ a.e.

\proclaim{Lemma 2.2}$T$ is elementary if and only if $T$ is
disjointness preserving and
$\sigma(X,X')-$con\-tinuous.\endproclaim

\demo{Proof}It is trivial that an elementary operator is
disjointness-preserving and continuous for $\sigma(X,X').$
For the converse we check that if $0\le f\le g\in X$ then
$0\le |Tf|\le |Tg|.$ It suffices by a density argument to
establish this when $f,g$ are both countably simple.  Pick a
maximal family of Borel sets $\{A_i:i\in I\}$ of finite
positive measure such that $|T(f\chi_{A_i})|\le
|T(g\chi_{A_i})|$ a.e.  This family is countable and so its
union $B$ is a Borel set.  If $A$ is a Borel set of positive
measure disjoint from $B$ then we may find a further set of
positive measure $A'\subset A$ so that $f\chi_{A'}=\alpha
g\chi_{A'}$ for some $0\le \alpha\le 1;$ thus
$|Tf_{\chi_{A'}}|\le |Tg_{\chi_{A'}}|$ contrary to our
maximality assumption.  Hence $\mu(\Omega\setminus B)=0.$
Now by Proposition 2.1, $Tf=\sum_{i\in I}T(f\chi_{A_i})$ and
$Tg=\sum_{i\in I}T(f\chi_{A_i})$ in $L_1(hd\mu)$ where $h$
is any strictly positive function in $X'.$ Thus $|Tf|\le
|Tg|.$ It follows that $T$ is regular and order-continuous
as an operator from $X$ into $L_1(hd\mu).$ As shown by Weis
[29] this means that $T$ is elementary.\bull\enddemo

\proclaim{Lemma 2.3}Suppose $\Omega$ is uncountable.  If $X$
is a K\"othe function space on $(\Omega,\mu)$ and $T:X\to X$
is an invertible elementary operator then $T^{-1}$ is
elementary and $T$ can be represented in the form
$Tf(s)=a(s)f(\sigma(s))$ where $a$ is a nonvanishing Borel
function and $\sigma:\Omega\to\Omega$ is an invertible Borel
map such that $\mu(B)>0$ if and only if
$\mu(\sigma^{-1}B)>0.$ \endproclaim

\demo{Proof}As noted above, $T$ can be represented in the
form $Tf=af\circ\sigma_0$ where $a$ is a Borel function and
$\sigma_0:\Omega\to\Omega$ is a Borel map.  Since $T$ is
onto it is clear that $a$ can vanish only on a set of
measure zero and so we may assume that it is nonvanishing.
Then for any $f,$ supp $Tf=\sigma_0^{-1}$(supp $f.)$ Thus
$T^{-1}$ is disjointness-preserving.  Now suppose $0\le
g_n\uparrow g$ a.e.; we will verify condition (2) of
Proposition 2.1 for $T^{-1}.$ We can suppose
$g_n=af_n\circ\sigma_0$ and $g=af\circ\sigma_0.$ Then
$T^{-1}g_n=f_n;$ we will show that, almost everywhere, we
have both $f_n(\omega)\to f(\omega)$ and $|f_n(\omega)|\le
|f(\omega)|$ for all $n.$ Once this is established then the
Dominated Convergence Theorem establishes Proposition 2.1
(2).  Suppose $E$ is any Borel set of finite measure such
that for every $\omega\in E$ we have
$\sup_n|f_n(\omega)|>|f_n(\omega)|$ or $f_n(\omega)$ does
not converge to $f(\omega).$ Then $\sigma_0^{-1}E$ is
contained in the set where $g_n(\omega)$ fails to converge
monotonically to $g(\omega)$ and so has measure zero.  This
implies that $T\chi_E=0$ (a.e.) and so $\mu(E)=0.$ Hence
$T^{-1}$ satisfies (2) of Proposition 2.1 and hence is
$\sigma(X,X')-$continuous.  It follows that $T^{-1}$ is
elementary and so can be represented in the form
$T^{-1}f=bf\circ\tau$ where $b$ is a nonvanishing Borel
function and $\tau:\Omega\to\Omega$ is a Borel map.  Thus
the identity map can be written in the form $f\to
a(b\circ\sigma_0)f\circ\tau\circ\sigma_0$ and so
$\tau\sigma_0s=s$ a.e.; similarly $\sigma_0\tau s=s$ a.e.

Let $E=\{s:\tau\sigma_0(s)=s\}.$ By Lusin's theorem there is
an increasing sequence of compact subsets $K_n$ of $E$ so
that $\sigma_0$ is continuous on $K_n$ and
$\mu(\Omega\setminus K)=0$ where $K=\cup K_n.$ Then
$\sigma_0$ is a Borel isomorphism of $K$ onto $\sigma_0(K)$
and both sets are $F_{\sigma}$'s.  Let $F$ be an uncountable
compact subset of $K$ of measure zero.  Then we define
$\sigma=\sigma_0$ on $K\setminus F$ and $\sigma=\rho$ on
$F\cup(\Omega\setminus K)$ where $\rho$ is any Borel
isomorphism between the two uncountable Borel sets
$F\cup(\Omega\setminus K)$ and
$\sigma_0(F)\cup(\Omega\setminus\sigma(K)).$ Then
$\sigma=\sigma_0$ a.e. and is a Borel automorphism.  We thus
can replace $\sigma_0$ by $\sigma$ and assume that $\sigma$
is a Borel automorphism.  Finally to show the measure
properties of $\sigma$ note that $\mu(B)=0$ if and only if
$T\chi_B=0$ a.e. i.e. if and only if
$\mu(\sigma^{-1}B)=0.$\bull\enddemo

\proclaim{Lemma 2.4}If $T:X\to X$ is an invertible
elementary operator then $T':X'\to X'$ is an elementary
operator.\endproclaim

\demo{Proof}We can represent $T$ in the form
$Tf=af\circ\sigma$ where $a$ is nonvanishing and $\sigma$ is
an invertible Borel map with $\mu(\sigma^{-1}B)=0$ if and
only if $\mu(B)=0.$ Let $w$ be the Radon-Nikodym derivative
of the $\sigma-$finite measure $\nu(B)=\mu(\sigma^{-1} B).$
Then for $f\in X,g\in X'$ we have $$ \align \int (T'g)fd\mu
&= \int gaf\circ\sigma d\mu\\ &= \int
(g\circ\sigma^{-1})(a\circ\sigma^{-1})fd\nu\\ &= \int
(g\circ\sigma^{-1})(a\circ\sigma^{-1})fwd\mu.  \endalign $$
Thus $T'g = a\circ\sigma^{-1}wg\circ\sigma^{-1}$ a.e. and
thus is elementary.

Of course if $X$ is order-continuous every operator $T:X\to
X$ is $\sigma(X,X')-$continuous.  However, for isometries we
can prove a similar result even without this assumption.

\proclaim{Proposition 2.5}Let $X$ be a K\"othe function
space with the Fatou property and suppose $T:X\to X$ is a
surjective isometry.  Then $T$ is $\sigma(X,X')-$continuous.
\endproclaim

\demo{Proof}We will use ideas developed in [9].  We recall
that the ball topology on $X$ is the weakest topology $b_X$
for which every closed ball (with any center and radius) is
closed.  Then $T:(X,b_X)\to (X,b_X)$ is continuous.  The
topology is not a Hausdorff topology, but ([9], Theorem 3.3)
its restriction to any absolutely convex Rosenthal set is
Hausdorff.  Here a set is a Rosenthal set if every sequence
contains a weakly Cauchy subsequence.

Now suppose $h$ is any strictly positive function in $X'.$

We show that if $0\le f_n\in X$ and $f_n\uparrow f$ a.e.
where $f\in X$ then $Tf_n$ converges to $Tf$ in
$L_1(hd\mu).$ In fact, setting $f_0=0,$ $\sum_{n\ge 1}
(f_n-f_{n-1})$ is weakly unconditionally Cauchy in $X$ and
so $\sum_{n\ge1} (Tf_n-Tf_{n-1})$ is weakly unconditionally
Cauchy in $X.$ Thus $\sum_{n\ge 1}(Tf_n-Tf_{n-1})$ converges
unconditionally to some $g$ in $L_1(hd\mu)$.

In particular $Tf_n$ converges in $L_1(h\mu)$ to $g.$ Since
$(Tf_n)$ is bounded in $X$ and $X$ has the Fatou property
(i.e.  $B_X$ is $L_0$-closed) it follows that $g\in X.$

Now consider the absolutely convex hull of $(Tf_n)$ together
with the points $g$ and $Tf$.  This is a Rosenthal set.
Since $b_X$ is weaker than the $L_0-$topology it follows
that $Tf_n$ converges to $g$ in $b_X.$ However $f_n$
converges to $f$ in $b_X$ and $Tf_n$ also converges to $Tf.$
We conclude that $Tf=g$ and so $(Tf_n)$ converges to $Tf$ in
$L_1(hd\mu).$

We now can conclude the argument by appealing to Proposition
2.1.\bull\enddemo

\demo{Remark}This result will only be needed to prove the
main result for nonseparable r.i. spaces.  The reader who is
only concerned with the separable case can observe that if
$X$ is separable it must be order-continous and then any
isometry $T$ is $\sigma(X,X')-$continuous.  Also its adjoint
$T':X'\to X'$ can be shown directly to be $\sigma(X',X'')$
continuous by identifying $X''$ as the sequential closure of
$X$ in $X^{**}.$\enddemo

\vskip2truecm

\subheading{3.  Flinn elements}

Let $X$ be a real Banach space and suppose $T:X\to X$ is a
linear operator.  We define $\Pi(X)$ to be the subset of
$X\times X^*$ of all $(x,x^*)$ such that
$\|x\|=\|x^*\|=x^*(x)=1.$ We recall that an operator $T:X\to
X$ is {\it numerically positive} (Rosenthal [26]) if
$x^*(Tx)\ge 0$ whenever $(x,x^*)\in\Pi(X).$ This is
equivalent to requiring the slightly weaker condition that
given $x$ with $\|x\|=1$ there exists $x^*$ so that
$(x,x^*)\in\Pi(X)$ and $x^*(Tx)\ge 0$ (see Lumer [20], [2]).
By results of Lumer [20] and Lumer and Phillips [23] (see
also [2]) it is equivalent to the requirement that
$\|\exp(-\alpha T)\|\le 1$ for $\alpha\ge 0.$ In the case
when $T$ is a projection it is easily seen that $T$ is
numerically positive if and only if $\|I-T\|=1.$

We next introduce an idea which is a real analogue of the
notion of hermitian elements [17].  Based on ideas of P.H.
Flinn [26] we say that $u\in X$ is a {\it Flinn element} if
there is a numerically positive projection $P:X\to [u].$ The
set of Flinn elements will be denoted $\Cal F(X).$ Note that
$0\in\Cal F(X)$ and that $u\in\Cal F(X)$ and $\alpha\in\bold
R$ imply $\alpha u\in \Cal F(X).$ If $0\neq u\in \Cal F(X)$
then there exists $f\in X^*$ so that $f\otimes u$ is
numerically positive projection onto $[u].$ We say then that
$(u,f)$ is a {\it Flinn pair.} Clearly $(u,f)$ is a Flinn
pair if and only if $f(u)=1$ and $f(x)x^*(u)\ge 0$ for
$(x,x^*)\in\Pi(X).$

\proclaim{Proposition 3.1}The set $\Cal F(X)$ is
closed.\endproclaim

\demo{Proof}Suppose $u_n\in\Cal F(X)$ and $\lim
\|u_n-u\|=0.$ It suffices to consider the case when
$\|u_n\|\neq 0$ and $\|u\|\neq 0.$ Then there exist $f_n\in
X^*$ so that $f_n\otimes u_n$ is a numerically positive
projection.  Thus $\|f_n\otimes u_n\|=\|f_n\|\|u_n\| \le 2.$
Thus $\|f_n\|\le 2\sup(1/\|u_n\|).$ By Alaoglu's theorem
$(f_n)$ has a weak$^*$-cluster point $f$ and clearly $(u,f)$
is a Flinn pair.\bull\enddemo

The next proposition is trivial, but we record it for future
use.

\proclaim{Proposition 3.2}Suppose $U:X\to Y$ is a surjective
isometry.  Then $U(\Cal F(X))= \Cal F(Y)$; furthermore if
$(u,f)$ is a Flinn pair then $(U(u),(U^*)^{-1}f)$ is a Flinn
pair.\endproclaim

The next theorem is due to Flinn (see [26], Theorem 1.1).

\proclaim{Theorem 3.3}Let $X$ be a Banach space and $\pi$ be
a contractive projection on $X$ with range $Y.$ Suppose
$(u,f)$ is a Flinn pair in $X$.  Suppose $f\notin
Y^{\perp}.$ Then $\pi(u)\in \Cal F(Y).$\endproclaim

\demo{Proof}Let $g$ be the restriction of $f$ to $Y.$ We may
assume $\pi(u)\neq 0.$ Then let $S=g\otimes \pi u$ be a rank
one operator on $Y.$ If $(y,y^*)\in \Pi(Y)$ then
$(y,y^*\circ \pi)\in \Pi(X).$ Now $f(y)y^*(\pi u)\ge 0$ and
so $S$ is numerically positive.  But $S^2=\beta S$ where
$\beta=g(\pi u).$ By considering $y=\pi u/\|\pi u\|$ and
choosing $y^*$ to norm $\pi u$ it is immediately clear that
$\beta\ge 0.$ If $\beta=0$ then $\exp(-\alpha S)=I-\alpha
S$; since by assumption $S$ is non-zero this contradicts
$\|\exp(-\alpha S)\|\le 1$ for all $\alpha\ge 0.$ Hence
$\beta>0$ and $(\pi u,\beta^{-1}g)$ is a Flinn
pair.\bull\enddemo

\vskip2truecm

\subheading{4.  Flinn elements in lattices}

Now suppose that $\Omega$ is a Polish space and that $\mu$
is a $\sigma-$finite Borel measure on $\Omega$.

\proclaim{Proposition 4.1} Let $X$ be an order-continuous
K\"othe function space on $\Omega$.  \newline (a) Suppose
that $(u,f)$ is a Flinn pair with $u\in X$ and $f\in
X'=X^*.$ Then $fu\ge 0$ a.e.\newline (b) Suppose $u\in\Cal
F(X)$.  Then there exists $f\ge 0$ such that $(|u|,f)$ is a
Flinn pair.\endproclaim

\demo{Proof}Let $A$ be a Borel subset of the set
$\{f>0\}\cap \{u<0\}$ of finite measure.  Suppose $\mu(A)>0$
and let $x=\chi_A/\|\chi_A\|.$ Pick $x^*$ so that
$(x,x^*)\in \Pi(X)$ and supp $x^*\subset A$.  Then $x^*\ge
0$ (a.e.) and $\int ux^*d\mu<0$ but $\int fx\,d\mu>0.$ This
contradiction shows that $\mu(A)=0$ and so the set
$\{f>0\}\cap \{u<0\}$ has measure zero.  Similar reasoning
shows that the set $\{f<0\}\cap \{u>0\}$ has measure zero.

(b) There is an isometry of $X$ onto $X$ which carries $u$
to $|u|$ so that $|u|\in\Cal F(X)$ by Proposition 3.2.  Now
suppose $(|u|,f)$ is a Flinn element.  Let $A=\{f<0\}$ and
consider the isometry $Ux =x-2\chi_Ax.$ Clearly by (a),
$U(|u|)=|u|$ and of course $(U^*)^{-1}f=|f|$ so that
$(|u|,|f|)$ is a Flinn pair.\bull\enddemo

\proclaim{Lemma 4.2} Suppose $\mu$ is nonatomic and suppose
$f,g\in L_1(\mu)$ with $\int |f|\,d\mu>0$ satisfy the
criterion that $$ \left(\int hf\,d\mu\right)\left(\int
hg\,d\mu\right)\ge 0$$ whenever $|h|=1$ a.e.  Then there is
a nonnegative constant $c$ so that $g=cf$ a.e.\endproclaim

\demo{Proof}Consider the subset $\Gamma$ of $\bold R^2$ of
all $(a,b)$ such that for some $h\in L_{\infty}(\mu)$ with
$|h|=1$ a.e. we have $\int hf\, d\mu=a$ and $\int
hg\,d\mu=b.$ Then it is an immediate consequence of
Liapunoff's theorem [27] that $\Gamma$ is closed and convex.
However $\Gamma=-\Gamma$ and the criterion is that $\Gamma$
is contained entirely in the union of the first and third
quadrants.  This trivially implies that $\Gamma$ is
contained in a line through the origin; the hypothesis on
$f$ implies this line is not the $y$-axis and so we deduce
the existence of $c\ge 0$ so that $\int hg\,d\mu=c\int
hf\,d\mu$ for all such $h$ and the lemma
follows.\bull\enddemo

We now establish the analogue of Theorem 6.5 of [17].

\proclaim{Theorem 4.3}Suppose $\mu$ is nonatomic and suppose
$X$ is an order-continuous K\"othe function space on
$(\Omega,\mu).$ Then $u\in X$ is a Flinn element if and only
if there is a nonnegative function $w\in L_0(\mu)$ with supp
$w=$ supp $u=B,$ so that:\newline (a) If $x\in X(B)$ then
$\|x\| =(\int |x|^2w\,d\mu)^{1/2}.$ \newline and\newline (b)
If $v\in X(\Omega\setminus B)$ and $x,y\in X(B)$ satisfy
$\|x\|=\|y\|$ then $\|v+x\|=\|v+y\|.$\endproclaim

\demo{Proof}Assume first that $0\neq u\in\Cal F(X)$.  We can
assume there exists $f\in X^*$ so that $(u,f)$ is a Flinn
pair.  Suppose first that $(x,x^*)\in \Pi(X).$ Then if
$|h|=1$ a.e. we also have $(hx,hx^*)\in \Pi(X)$ and so
$(\int uhx^*\,d\mu)(\int fhx\,d\mu)\ge 0.$ By Lemma 4.2,
there is a constant $k_x>0$ so that $ux^*=k_xfx$ almost
everywhere.  It follows immediately that we must have
$f\chi_{\Omega\setminus B}=0$ almost everywhere.  Thus we
can define a function $w$ by $w=f/u$ on $B$ and $w=0$
otherwise.  Then if $(x,x^*)\in \Pi(X)$ we have
$x^*\chi_B=k_xwx\chi_B.$

Next let us suppose that $e_1,e_2\in X(B)$ satisfy the
conditions $\int e_1^2w\, d\mu =\int e_2^2w\, d\mu=1$ and
$\int e_1e_2w\,d\mu=0.$ Consider the function
$F(\varphi)=e_1\cos\varphi+e_2\sin\varphi$ for $0\le
\varphi\le 2\pi.$ Suppose $v\in X(\Omega\setminus B)$ and
consider the function $H(\varphi)=\|v+F(\varphi)\|.$

We note that the function $H$ is Lipschitz on $[0,2\pi].$ We
will show that $H'(\varphi)=0$ a.e. and deduce that $H$ is
constant.

Let us suppose that $\theta$ is a point of differentiability
of $H$.  Let $g\in X^*$ be a norming function for
$v+F(\theta).$ Then $H(\varphi)-\langle
v+F(\varphi),g\rangle$ has a minimum at $\varphi=\theta$ and
so we can deduce tht $H'(\theta)=\langle
F'(\theta),g\rangle.$

Since $g$ norms $v+F(\theta)$ we conclude that
$g\chi_B=cwF(\theta)$ for some nonnegative constant $c.$
Thus $$ \align \langle F'(\theta),g\rangle &= c\int w
F(\theta)F'(\theta)d\mu\\ &= c\cos2\theta\int we_1e_2
d\mu-\frac{c}2\sin2\theta \int w(e_1^2-e_2^2)d\mu\\ &=0.
\endalign $$ Thus $H$ is constant as promised.

It follows immediately that if $x,y\in X(B)$ satisfy $\int
x^2w\,d\mu=\int y^2w\,d\mu=1$ then $\|v+x\|=\|v+y\|$; simply
determine $e_2$ so that $\int e_2xw\,d\mu=0$, $\int
e_2^2wd\mu=1$ and $y=x\cos\varphi +e_2\sin\varphi.$

Taking the special case $v=0$ this leads easily to (a).  (b)
is then the general case.

The converse is easy.  First note that (b) easily implies
that if $x,y\in X(B)$ with $\|x\|\le \|y\|$ then for $v\in
X(\Omega\setminus B)$ we have $\|v+x\|\le \|v+y\|.$ Suppose
$\|u\|=1$ and (a) and (b) hold.  We show that the pair
$(u,uw)$ is Flinn.  Clearly $\langle u,uw\rangle=1.$ Suppose
$x\in X$; then $I-uw\otimes u(x)= x\chi_{\Omega\setminus
B}+y$ where $\|y\|\le \|x\chi_B\|$ and so $\|I-uw\otimes
u\|\le 1.\bull$\enddemo

We now apply this theorem to the case when $X$ is a
separable r.i. space on $[0,1].$

\proclaim{Theorem 4.4}Suppose $X$ is a separable r.i. space
on $[0,1]$.  If $\Cal F(X)\neq \{0\}$ then $X=L_2[0,1].$
\endproclaim

\demo{Proof}If $\Cal F(X)\neq 0$ then Theorem 4.3 shows that
there is a Borel set $B\subset [0,1]$ of positive measure
and a weight function $w\in L_0$ such that if supp $f\subset
B$ then $\|f\|_X=(\int |f|^2w d\lambda)^{1/2}.$ Further if
$g\chi_B=0$ and $f_1,f_2\in X(B)$ then
$\|g+f_1\|_X=\|g+f_2\|_X.$ It follows immediately from
re-arrangement invariance that $w$ is constant and we obtain
the existence of $c,\delta>0$ so that if $\lambda(\text{supp
}f)\le \delta$ then $\|f\|_X =c\|f\|_2.$

Now pick an integer $N$ so that $1/N<\delta.$ It follows
easily from condition (b) of the previous theorem that there
is a constant $a>0$ so that if $f_1,f_2,\ldots,f_N$ are
disjoint functions satisfying $\lambda(\text{supp }f_k)\le
1/N$ and $\|f_k\|_2=1$ then $\|f\|_X=a.$ Consider then any
simple function $f$ and write $f=f_1+\cdots+f_N$ where $f_k$
are identically distributed and disjointly supported.  Then
$\|f\|_X=a\|f_1\|_2=aN^{-1/2}\|f\|_2.$ By considering
$\chi_{[0,1]}$ it is clear that $aN^{-1/2}=1$ and the
theorem follows easily.\bull\enddemo

\vskip2truecm

\subheading{5.  Flinn elements of finite-dimensional r.i.
spaces}

Suppose $N$ is a natural number.  Let
$e^N_i=\chi_{((i-1)2^{-N},i2^{-N}]}$ for $1\le i\le 2^N.$
Let $X_N=[e^N_i:1\le i\le 2^N].$ We denote the averaging
projection (conditional expectation operator) of $X$ onto
$X_N$ by $\Cal E_N.$ Notice that $X_N^*$ can be identified
naturally with $X'_N.$ We will also let $X_N^{-}=[e^N_i:1\le
i\le 2^N-1].$

\proclaim{Lemma 5.1}Suppose $X$ is an r.i. space on $[0,1]$
so that $X\neq L_2.$ Then there exists $N\in\bold N$ so that
$\sum_{i<2^N}e^N_i \notin \Cal F(X_N^{-}).$\endproclaim

\demo{Proof}Suppose for every $n\in \bold N$ we have
$\chi_n=\sum_{i<2^n}e^n_i=\chi_{[0,1-2^{-n}]}\in\Cal
F(X_n^{-}).$ Since $(X_n^-)^*$ can be identified with
$(X'_n)^{-}$ there exists $f_n=\sum_{i<2^n}a_{ni}e^n_i$ so
that $(\chi_n,f_n)$ is a Flinn pair for $X_n^{-}$ i.e.
$\int f_n\,d\lambda=1$ and $\|I-f_n\otimes \chi_n\|=1.$ Then
for every permutation $\sigma$ of ${1,2,\ldots,2^n-1}$ we
have that $(\chi_n,f_n^{\sigma})$ is a Flinn pair where $
f_n^{\sigma}=\sum_{i<2^n}a_{n\sigma(i)}e_i^n.$ By averaging
we conclude that $(\chi_n,(1-2^{-n})^{-1}\chi_n)$ is a Flinn
pair.

Now suppose $x\in X$.  We conclude that $$ \|\Cal
E_n(x\chi_n) - (1-2^{-n})^{-1}(\int_0^{1-2^{-n}}
x(t)\,dt)\chi_n\|_X \le \|\Cal E_n(x\chi_n)\|_X.$$

Letting $n\to \infty$ we obtain (by the Fatou property of
the norm when $X$ is not separable) that $$ \|x - (\int_0^1
x(t)\,dt)\chi_{[0,1]}\|_X \le \|x\|_X$$ and so
$(\chi_{[0,1]},\chi_{[0,1]})$ is a Flinn pair in $X\times
X'.$ Now if $X$ is separable (i.e. order-continuous) Theorem
4.4 gives the conclusion that $X$ is isometric to $L_2$.  If
not we consider $X_0$, the closure of the simple functions
in $X$; it is immediate that $\chi_{[0,1]}$ is Flinn in
$X_0$ and so if $X_0$ is separable, we can again apply
Theorem 4.4 to get the conclusion that $X$ is isometric to
$L_2.$ There remains one case, when $X_0$ is not
order-continuous and so ([19]) $X_0=L_{\infty}[0,1]$ up to
renorming.  But then we conclude that $\chi_{[0,1]}$ is
Flinn in $X'$ which is $L_1$ up to renorming and get a
contradiction.  \bull\enddemo

We now need to introduce a technical definition.  We will
say that an r.i. space $X$ has property $(P)$ if for every
$t>0,$ $$ \|e^1_1\|_X < \|e^1_1+te^1_2\|_X.$$ We say that
$X$ has property $(P')$ if $X'$ has property (P).

\proclaim{Lemma 5.2}Any r.i. space $X$ has at least one of
the properties $(P)$ or $(P').$ \endproclaim

\demo{Proof}Assume $X$ fails both $(P)$ and $(P').$ Then for
small enough $\eta>0$ we have $\|e^1_1+\eta
e^1_2\|_X=\|e^1_1\|_X$ and $\|e^1_1+\eta e^1_2\|_{X'}
=\|e^1_1\|_{X'}.$ But then $$ \align \frac12(1+\eta^2)
&=\int (e_1^1+\eta e_2^1)^2 d\lambda\\ &\le
\|e^1_1\|_X\|e^1_1\|_{X'}\\ &=\frac12.  \endalign $$ This
contradiction establishes the lemma.\bull\enddemo

\demo{Remark}If $X$ is strictly convex then it has property
$(P).$

\proclaim{Lemma 5.3}Assume $X$ has property $(P).$ Suppose
$(e^N_j,u)$ is a Flinn pair in $X_N\times X'_N.$ Then
$u=2^Ne^N_j.$\endproclaim

\demo{Proof}It suffices to consider the case $j=1$.  We can
write $u=2^Ne^N_1 + \sum_{j>1}a_je^N_j.$ By using
Proposition 3.2 it follows that $(e^N_1,|u|)$ is also a
Flinn pair.  Then by an averaging procedure as in the
preceding Lemma 5.1 we can show that $(e^N_1,v)$ is a Flinn
pair where $v=2^Ne^N_1 + \eta\sum_{j\ge 2}e^N_j$ where
$(2^N-1)\eta=\sum_{j\ge 2}|a_j|.$ We now project by $\Cal
E_1$ onto $X_1.$ By Theorem 3.3, $(\Cal E_1e^N_1,w)$ is a
Flinn pair where $w$ is a multiple of $\Cal E_1v.$ Thus
$(e^1_1,2(e^1_1+\tau e^1_2))$ is a Flinn pair for some
$\tau>0.$

Now consider $g=\frac12\tau e^1_1 - e^1_2\in X_1$ and
suppose this is normed by $h= \alpha e^1_1-\beta e^1_2\in
X'_1$, where $\alpha,\beta\ge 0.$ Thus $\|h\|_{X'}=1$ and
$\frac14\tau\alpha +\frac12\beta= \|g\|_X.$ Now $\int
2g(e^1_1+\tau e^1_2)d\lambda = -\frac12\tau<0,$ and hence
$\int he^1_1d\lambda\le 0$ i.e.  $\alpha\le 0.$ hence
$\alpha=0$ and so $h=-\|e^1_2\|_{X'}^{-1}e^1_2$ and $\|g\|_X
= \|e^1_2\|_X$ which contradicts property
$(P).$\bull\enddemo

\proclaim{Lemma 5.4}Suppose $N$ is a natural number and
$N=lm.$ Suppose $d_1\ge d_2\ge\cdots\ge d_N\ge 0.$ Then
there is a permutation $\sigma$ of $\{1,2,\ldots,N\}$ so
that if $$b_j= \sum_{i=1}^m d_{\sigma((j-1)m+i)}$$ for $1\le
j\le l$ then $\max_{i,j}|b_i-b_j|\le d_1.$\endproclaim

\demo{Proof}The construction is inductive.  We will define
$\sigma((j-1)m+k)$ in blocks for $k=1,2,\ldots,m.$ For $k=1$
we define $\sigma((j-1)m+1)=j.$ Now suppose we have
completed the construction up to $k-1<m$.  We calculate
$b^{k-1}_j=\sum_{i=1}^{k-1}d_{\sigma((j-1)m+i}.$ We then
define $\sigma((j-1)m+k)\in [(k-1)l+1,kl]$ in such a way
that $b^{k-1}_i < b^{k-1}_j$ implies $\sigma((i-1)m+k)<
\sigma((j-1)m+k).$ This describes the construction of
$\sigma.$

Now by induction we have $\max_{i,j}|b^k_i-b^k_j|\le d_1.$
For $k=1$ this is obvious.  Suppose we have the result for
$k-1.$ Suppose $\sigma((i-1)m+k)\le \sigma((j-1)m+k).$ Then
$b^{k-1}_i\le b^{k-1}_j;$ furthermore $b^k_i=b^{k-1}_i+x$
and $b^k_j=b^{k-1}_j=y$ where $0\le y\le x$.  Hence
$|b^k_i-b^k_j| \le \max(x-y,b^{k-1}_j-b^{k-1}_i)\le d_1.$

Now if we let $k=l,$ the lemma is proved.  \bull\enddemo

\proclaim{Proposition 5.5}Suppose $X$ is an r.i. space on
$[0,1]$ with property $(P')$ and such that $X\neq L_2.$ Then
for any $0<p<\infty$ there is a constant $A_p=A_p(X)$ so
that for every $n\in\bold N$ and every
$u=\sum_{i=1}^{2^n}a_ie^n_i\in\Cal F(X_n)$ we have $$
\left(\sum_{i=1}^{2^n}|a_i|^p\right)^{1/p} \le A_p\max_{1\le
i\le 2^n}|a_i|.$$ \endproclaim

\demo{Proof}We start with the simple observation that if
$\sum a_ie^n_i$ is Flinn then so is $\sum |a_i|e^n_i$ and so
it suffices to consider only the case when $u\ge 0.$
Similarly we are free to permute the $(a_i)$.  We therefore
consider the case when $a_1\ge a_2\ge \cdots \ge a_{2^n}\ge
0.$

Now according to Lemma 5.1 there exists $m$ so that
$\sum_{i<2^m}e^m_i\notin \Cal F(X_m^{-}).$ In fact by
Proposition 3.1, this means that there exists $\delta>0$ so
that if $w\in\Cal F(X_m^{-}))$ then
$\|w-(2^m-1)^{-1}\sum_{i<2^m}e^m_i\|_{\infty}\ge \delta/2.$
This implies that if $w=\sum_{i<2^m}b_ie^m_i$ and $\sum
b_i=1$ then $\max_{i,j}|b_i-b_j|\ge \delta.$

Now let us suppose $n> m$ and that
$u=\sum_{j=1}^{2^n}a_je^n_j$ is Flinn in $X_n$ where $a_1\ge
a_2\ge \cdots \ge a_{2^n}\ge 0.$ Let us set $S_k=\sum_{j\le
k}a_j$ for $1\le k\le 2^n.$ Let $S_0=0$ and $S=S_{2^n-2^m}.$

Fix $1\le k\le 2^{n-m}.$ We consider a permutation $\sigma$
of $\{1,2,\ldots,2^n\}$ so that
$\sigma\{2^n-2^{n-m}+1,\ldots,2^n\}=\{i:i<k\}\cup\{i:i\ge
2^n-2^{n-m}+k\}$ and such that if $$ b_j
=\sum_{i=1}^{2^{n-m}}a_{\sigma((j-1)2^{n-m}+i)}$$ for $1\le
j\le 2^m-1$ then $\max|b_i-b_j| \le a_k.$ Such a permutation
exists by Lemma 5.4.

Now we argue that if $v=\sum_{j=1}^{2^n}a_{\sigma(j)}e^n_j$
then $v\in \Cal F(X_n)$ and so $\Cal E_m(v)\in\Cal F(X_m)$.
To see this observe that there exists $g\in X'_n$ with $g\ge
0$ so that $(v,g)$ is a Flinn pair by Proposition 4.1;
clearly $\Cal E_m(g)\neq 0$ and so by Theorem 3.3, $\Cal
E_m(v)\in \Cal F(X_m).$ Thus $w=\sum_{i\le 2^m}b_ie^m_i \in
\Cal F(X_m).$

Next we claim that $w_0=\sum_{i<2^m}b_ie^m_i\in \Cal
F(X_m^{-}).$ If $w_0=0$ this is trivial.  If not, select
$h\ge 0$ in $X'_m$ so that $(w,h)$ is a Flinn pair.  If
$h=\sum_{i\le 2^m}c_ie^m_i$ we argue that there exists
$i<2^m$ so that $c_i> 0.$ For, if not, $h$ is a multiple of
$e^m_{2^m}$ and by Lemma 5.3, since $X'$ has $(P)$, we get
that $b_i=0$ for $i<2^m,$ i.e.  $w_0=0.$ Now we can apply
Theorem 3.3 to deduce that $w_0\in\Cal F(X_m^{-}).$

Recalling the original choice of $\delta$ this implies that:
$$\max_{i,j<2^m}|b_i-b_j|\ge \delta\sum_{j=1}^{2^m-1}b_j.$$

In view of the selection of $\sigma$ we have $$ a_k \ge
\delta(S_{2^n-2^{n-m}+k-1}-S_{k-1})\ge \delta(S-S_{k-1})$$
and this holds for $1\le k\le 2^{n-m}.$ For convenience, let
us put $\alpha=1-\delta.$ Then, for $1\le k\le 2^{n-m}$ we
have $$ (S-S_k) \le \alpha (S-S_{k-1}).$$ By induction, we
have $$ (S-S_k)\le \alpha^kS$$ for $1\le k\le 2^{n-m}.$ This
gives an estimate on $a_k$, i.e.  $$ a_k \le S-S_{k-1} \le
\alpha^{k-1}S\le \delta^{-1}\alpha^{k-1}a_1,$$ for $1\le
k\le 2^{n-m}.$

If $0<p<\infty$, this implies that $$ \align
\sum_{i=1}^{2^n}a_i^p &\le 2^m\sum_{i=1}^{2^{n-m}}a_i^p\\
&\le 2^ma_1^p\delta^{-p}(1-\alpha^p)^{-1} \\ &= a_1^pB_p^p,
\endalign $$ say.  This estimate holds if $n>m.$ If we take
$A_p=\max(2^{m/p},B_p)$ we obtain the Proposition as stated.
\bull\enddemo

\newpage

\subheading{6.  Isometries on r.i. spaces}

\proclaim{Theorem 6.1}Let $X$ be an r.i. space on $[0,1]$
with $X\neq L_2.$ Suppose $X$ has property $(P).$ Then for
any $0<p\le 1$ there is a constant $C_p=C_p(X)$ with the
following property.  Suppose $Y$ is any K\"othe function
space on some Polish space $(\Omega,\mu),$ for which $Y'$ is
norming.  Suppose $T:X\to Y$ is an isometric isomorphism of
$X$ onto $Y.$ Then $$
\sup_n\|(\sum_{i=1}^{2^n}|Te^n_i|^p)^{1/p}\|_Y \le C_p.$$
\endproclaim

\demo{Remark} The sequence
$(\sum_{i=1}^{2^n}|Te^n_i|^p)^{1/p}$ is increasing.  If $Y$
has the Fatou property it will follow that
$\sup_n(\sum_{i=1}^{2^n}|Te_i^n|^p)^{1/p} \in Y.$ \enddemo

\demo{Proof} We note first that by Proposition 2.5, $T^{-1}$
is $\sigma(X,X')-$continuous and so has an adjoint
$S=(T^{-1})':X'\to X'.$ We define $f^n_i=Te^n_i$ and
$g^n_i=Se^n_i$.  Suppose $(x,x^*)\in \Pi(X_n)$ where $x=\sum
a_ie^n_i$ and $x^*=\sum a_i^*e_i^n.$ Then $(Tx,Sx^*)\in
\Pi(Y)$ and this implies that $$
(\sum_{i=1}^{2^n}a_if^n_i(\omega))(\sum_{i=1}^na_i^*g^n_i(\omega))\ge
0\tag *$$ for $\mu-$a.e.  $\omega\in\Omega.$

Using the fact that $\Pi(X_n)$ is separable it follows that
there is a set of measure zero $\Omega^n_0$ so that if
$\omega\notin \Omega^n_0,$ (*) holds for every $(x,x^*)\in
\Pi(X_n).$ Let $\Omega_0=\cup_{n\ge 1}\Omega_0^n.$

Now define
$F_n(\omega)=\sum_{i=1}^{2^n}f^n_i(\omega)e^n_i\in X'_n$ and
$G_n(\omega)=\sum_{i=1}^{2^n}g^n_i(\omega)e^n_i\in X_n.$ The
above remarks show the operator $G_n(\omega)\otimes
F_n(\omega)$ is numerically positive on $X'_n$ if
$\omega\notin\Omega_0.$

Now let $B_n=\{\omega:G_n(\omega)=0\}.$ Clearly $(B_n)$ is a
descending sequence of Borel sets.  Let $B=\cap B_n.$ If
$\mu(B)>0$ then there exists a nonzero $h\in Y$ supported on
$B$ and $\langle h, Sx'\rangle=0$ for every $x'\in X'$.
Thus $T^{-1}h=0,$ which is absurd.

Let $D_n=\Omega\setminus(\Omega_0\cup B_n).$ If $\omega\in
D_n$ then $G_n(\omega)\neq 0$ and so it follows that
$F_n(\omega)\in\Cal F(X'_n).$ We recall that $X$ has
property $(P)$ and so $X'$ has property $(P').$ Hence
letting $A_p=A_p(X')$ be the constant from Proposition 5.5
$$ (\sum_{i=1}^{2^n}|f^n_i(\omega)|^p)^{1/p} \le
A_p\max_{1\le i\le n}|f^n_i(\omega)|.$$

Hence $$ \align
\|\chi_{D_n}(\sum_{i=1}^{2^n}|f^n_i|^p)^{1/p}\|_Y &\le
A_p\|\max_{1\le i\le 2^n}|f^n_i|\|_X\\ &\le A_p
\|(\sum_{i=1}^{2^n}|f^n_i|^2)^{1/2}\|_X\\ &\le K_GA_p
\|(\sum_{i=1}^{2^n}|e^n_i|^2)^{1/2}\|_X \\ &=K_GA_p
\endalign $$ by Krivine's theorem ([19] 1.f.14, p.93.)  Now
the sequence $\chi_{D_n}(\sum_{i=1}^{2^n}|f^n_i|^p)^{1/p}$
is increasing, as $0<p\le 1.$ If $g\ge 0$ and $\|g\|_{Y'}\le
1,$ we have $$ \int_{D_n}
g(\sum_{i=1}^{2^n}|f^n_i|^p)^{1/p}d\mu \le K_GA_p$$ and so
$$ \int_{\Omega}
g(\sup_n(\sum_{i=1}^{2^n}|f^n_i|^p)^{1/p})d\mu\le K_GA_p.$$
We now quickly obtain the Theorem since $Y'$ is norming.
\bull\enddemo

Let $\Cal M=\Cal M[0,1]=C[0,1]^*$ denote the space of
regular Borel measures on $[0,1]$.  If $0<p\le 1$ and
$\mu\in\Cal M$ we define the $p$-variation of $\mu$ by $$
\|\mu\|_p = \sup\{(\sum_{k=1}^n|\mu(B_k)|^p)^{1/p}:\
n\in\bold N,\ B_1,\ldots,B_n\in\Cal B \text{ disjoint}\}.$$
If $p=1$ this reduces to the usual variation norm.  For
$p<1$ it is easily seen that $\|\mu\|_p<\infty$ if and only
if $\mu=\sum_{n=1}^{\infty}a_n\delta(t_n)$ for some sequence
of distinct elements $(t_n)$ in $[0,1]$ and $a_n\in\bold R$
such that $\sum |a_n|^p=\|\mu\|_p^p$ (see [14], [24]).  The
following lemma is standard and we omit the proof.

\proclaim{Lemma 6.2}For $\mu\in\Cal M$ we have $$\|\mu\|_p
=\sup_n(\sum_{k=1}^{2^n}|\mu(D(n,k))|^p)^{1/p}$$ where
$D(n,1)=[0,2^{-n}]$ and $D(n,k)=((k-1)2^{-n},k2^{-n}]$ for
$2\le k\le 2^n.$\endproclaim

We now use the machinery developed in [15].  Suppose $X$ is
an r.i. space.  Let $T:X \to L_0[0,1]$ be a continuous
linear operator.  We say that $T$ is {\it controllable} if
there exists $h\in L_0$ so that $|Tx| \le h$ a.e. when
$\|x\|_{\infty}\le 1.$ $T$ is said to be {\it
measure-continuous} if it satisfies the criterion that that
$|Tx_n|$ converges to zero in $L_0$ whenever
$\sup\|x_n\|_{\infty}\le 1$ and $|x_n|$ converges to zero in
$L_0.$ Thus it follows from Theorem 3.1 of [15] (cf.
Sourour [28]) that $T$ is controllable and
measure-continuous if and only if there is a weak$^*$-Borel
map $s\to\nu^T_s$ from $[0,1]$ into $\Cal M$ satisfying
$|\nu^T_s|(B)=0$ almost everywhere when $B$ has measure
zero, and such that we have for any $x\in L_{\infty}$ $$
Tx(s) = \int x(t)\,d\nu^T_s(t).$$ The map $s\to \nu^T_s$ is
called the {\it representing kernel} or {\it representing
random measure} for $T$ and it is unique up to sets of
measure zero.

We further remark that if $\|\nu^T_s\|_p<\infty$ a.e. for
some $p<1$ then $\nu^T_s$ is purely atomic for almost every
$s$ and so (cf.  [14], [29] Theorem 4.1) there is a sequence
of Borel maps $\sigma_n:[0,1]\to[0,1]$ and Borel functions
$a_n$ on $[0,1]$ so that $|a_n(s)|\ge |a_{n+1}(s)|$ a.e. for
every $n$ and $\sigma_m(s)\neq \sigma_n(s)$ whenever $m\neq
n$ and $s\in [0,1]$ and for which $$ \nu^T_s=
\sum_{n=1}^{\infty}a_n(s)\delta(\sigma_n(s)).$$

We can now summarize our conclusions, restricting attention
to surjective isometries on $X.$

\proclaim{Proposition 6.3}Let $X$ be an r.i. space on
$[0,1]$ with property $(P)$, and such that $X\neq L_2$
(isometrically).  Then for any $0<p\le 1$ there is a
constant $C_p$ depending only on $X$ such that the following
holds.  Suppose $T:X\to X$ is a surjective isometry.  Then
$T$ is controllable and further its representing kernel
$\nu^T_s$ satisfies $$ \int_0^1 \|\nu^T_s\|_p^pds \le
C_p^p.$$ \endproclaim

\demo{Proof}This is an almost immediate consequence of
Theorem 6.1.  We use the same notation.  If
$F=\sup_n(\sum_{k=1}^{2^n}|Te^n_k|)$ then we have
$\|F\|_{X''} \le C_1$ and so $\|F\|_1\le C_1$.  It is easy
to deduce that if $\|x\|_{\infty}\le 1$ then $|Tx|\le F.$ To
conclude the argument we will need that $T$ is
measure-continuous.  This is immediate if $X$ is not equal
to $L_{\infty}$ with some equivalent renorming, since in
this case $\|x_n\|_{\infty}\le 1$ and $|x_n|\to 0$ in
measure imply that $\|x_n\|_X\to 0.$ In the exceptional case
we use Propositions 2.1 and 2.5 to deduce that $T$ is
measure-continuous.  We conclude that in every case $T$ has
a representing random measure $\nu^T_s.$

Now if $0<p\le 1$, then by Lemma 6.2 $$ \|\nu^T_s\|_p
=\sup_n (\sum_{k=1}^{2^n}|Te^n_k(s)|^p)^{1/p}$$ almost
everywhere.  Hence $$ (\int \|\nu^T_s\|^p_pds)^{1/p} \le
\|\|\nu^T_s\|_p\|_{X''} \le C_p.\bull$$ \enddemo

\proclaim{Theorem 6.4}Let $X$ be an r.i. space on $[0,1]$
which is not isometrically equal to $L_2[0,1],$ and let
$T:X\to X$ be a surjective isometry.  Then there exists a
Borel function $a$ on $[0,1]$ with $|a|>0$ and an invertible
Borel map $\sigma:[0,1]\to [0,1]$ such that
$\lambda(\sigma^{-1}(B))>0$ if and only if $\lambda(B)>0$
for $B\in\Cal B$ and so that $Tx(s) = a(s)x(\sigma(s))$ a.e.
for every $x\in X.$ \endproclaim

\demo{Proof} We start by assuming that $X$ has property
$(P).$ According to the previous proposition every
surjective isometry $T$ is controllable and further for
every $0<p\le 1$ there is a constant $C_p$ depending only on
$X$ so that $$ (\int_0^1\|\nu_s^T\|^p_p ds)^{1/p} \le C_p.$$
Let us define $K_p$ to be the least such constant i.e.
$$K_p=\sup\{\|(\|\nu_s^T\|_p)\|_p:\text{ T is a surjective
isometry}\}.$$

Suppose $T$ is any fixed isometry.  We can represent
$\nu^T_s =\sum_{n=1}^{\infty}a_n(s)\delta(\sigma_n(s))$
where $a_n$ is a sequence of Borel functions, and
$\sigma_n:[0,1]\to [0,1]$ is a sequence of Borel maps
satisfying $\sigma_m(s)\neq \sigma_n(s)$ whenever $m\neq n$
and $0\le s\le 1.$ In this representation we can assume that
$\sigma_i(s)\neq 0$ for all $i,s$ since the measure of set
where $\nu^T_s(\{0\})\neq 0$ is clearly zero; thus we could
simply redefine $\sigma_i$ to avoid $0$ without changing the
kernel except on a set of measure zero.  It follows that $$
\|\nu^T_s\|_p^p =\sum_{n=1}^{\infty}|a_n(s)|^p.$$ We define
the function $H_p(s)=\sum |a_n(s)|^p.$

>From now on we will fix $0<p\le 1.$ Let $M_N(s)$ be the
greatest index such that $\sigma_1(s),\ldots,\sigma_M(s)$
belong to distinct dyadic intervals $D(N,k).$ Then
$M_N(s)\to\infty$ for all $s$ and it follows easily that
given $\epsilon>0$ we can find $M,N$ and a Borel subset $E$
of $[0,1]$ with $\lambda(E)>1-\epsilon$ and such that
$M_N(s)\ge M$ for $s\in E,$ and $$ \align
\int_{[0,1]\setminus E}H_pdt &<\epsilon\\
\int_E\sum_{n=M+1}^{\infty}|a_n(s)|^p ds&<\epsilon.
\endalign $$

For notational convenience we will set $P=2^N.$ Let us
identify the circle group $\bold T$ with $\bold R/\bold
Z=[0,1)$ in the natural way.  For $\theta\in [0,1)^{P}$ we
define a measure preserving Borel automorphism
$\gamma=\gamma(\theta_1,\ldots,\theta_P)$ given by
$\gamma(0)=0$ and then $$
\gamma(s)=s+(\theta_k-\rho)2^{-N}$$ for $(k-1)2^{-N}<s\le
k2^{-N}$ where $\rho=1$ if $2^Ns+\theta_k>k$ and $\rho=0$
otherwise.  Thus $\gamma$ leaves each $D(N,k)$ invariant.
The set of all such $\gamma$ is a group of automorphisms
$\Gamma$ which we endow with the structure of the
topological group $\bold T^P=[0,1)^P.$ We denote Haar
measure on $\Gamma$ by $d\gamma.$ For each $k$ let
$\Gamma_k$ be the subgroup of all $\gamma(\theta)$ for which
$\theta_i=0$ when $i\neq k.$ Thus
$\Gamma=\Gamma_1\ldots\Gamma_P.$

We also let the finite permutation group $\Pi_P$ act on
$[0,1]$ by considering a permutation $\pi$ as inducing an
automorphism also denoted $\pi$ by $\pi(0)=0$ and then
$\pi(s)=\pi(k)-k+s$ for $(k-1)2^{-N}<s\le k2^{-N}.$ We again
denote normalized Haar measure on $\Pi_P$ by $d\pi.$ Finally
note that the set $\Gamma\Pi_P=\Cal T$ also forms a compact
group when we endow this with the product topology and Haar
measure $d\tau=d\gamma\,d\pi$ when $\tau=\gamma\pi.$

We now wish to consider the isometries $V_{\tau}:X\to X$ for
$\tau\in\Cal T$ defined by $V_{\tau}x=x\circ\tau.$ For every
$\tau\in\Cal T$ the operator $S(\tau)=TV_{\tau}T$ is a
surjective isometry and so has an abstract kernel
$\nu^{S(\tau)}_s.$ \enddemo

\proclaim{Lemma 6.5}For almost every $\tau\in \Cal T$ we
have that $$ \align
\int_0^1\sum_{n=1}^{\infty}\sum_{j=1}^{\infty}|a_j(s)||a_n(\tau\sigma_js)|
ds &<\infty \tag 1\\
\sum_{n=1}^{\infty}\sum_{j=1}^{\infty}a_j(s)a_n(\tau\sigma_js)\delta(
\sigma_n\tau \sigma_is) &= \nu^{S(\tau)}_s\text{ a.e.}\tag 2
\endalign $$ \endproclaim

\demo{Proof}Let us prove (1).  Note that for any fixed $s$
and $(n,j)$ we have $$ \int_{\Cal
T}|a_n(\tau\sigma_js)|d\tau = \int_0^1 |a_n(t)|dt$$ and so
it follows that $$ \int_0^1\int_{\Cal T}
\sum_{n=1}^{\infty}\sum_{j=1}^{\infty}|a_j(s)||a_n(\tau\sigma_js)|
d\tau\,ds <\infty.$$ This proves the first assertion.  Note,
in particular, if (1) holds,
$$\sum_{n=1}^{\infty}\sum_{j=1}^{\infty}|a_j(s)||a_n(\tau\sigma_js)|<\infty$$
for almost every $s.$

To obtain (2) let us suppose that $\tau$ is such that (1)
holds.  Suppose $\|x\|_{\infty}\le 1.$ Then $V_{\tau}Tx$ may
not be bounded but there is an increasing sequence $F_m$ of
Borel subsets of $[0,1]$ with $\cup F_m=[0,1],$ so that
$\chi_{F_m}V_{\tau}Tx$ is bounded.  Thus $$ \align
T(\chi_{F_m}V_{\tau}Tx)(s) &=
\sum_{j=1}^{\infty}a_j(s)\chi_{\sigma_j^{-1}F_m}(s)Tx(\tau\sigma_j
s) \\ &=
\sum_{j=1}^{\infty}\left(\sum_{n=1}^{\infty}a_j(s)\chi_{\sigma_j^{-1}F_m}(s)
a_n(\tau\sigma_js)x(\sigma_n\tau\sigma_js) \right).
\endalign $$

Now for almost every $s$ since the double series absolutely
converges we may obtain $$
\lim_{m\to\infty}T(\chi_{F_m}V_{\tau}Tx)(s)=\sum_{n=1}^{\infty}\sum_{
j=1}^{\infty}a_j(s)a_n(\tau\sigma_js)x(\sigma_n\tau\sigma_js).$$

If $X$ is order-continuous of $X$ the left hand side is
simply $(S(\tau)x)(s)$ a.e.  Thus by the uniqueness of the
representing random measure we obtain (2).  For the general
case we use Propositions 2.1 and 2.5 to give the same
conclusion.  It follows that $$TV_{\tau}Tx(s) =
\sum_{n=1}^{\infty}\sum_{j=1}^{\infty}a_j(s)a_n(\tau\sigma_js)x(
\sigma_n\tau\sigma_js)$$ for almost every $s.$ Again the
uniqueness of the representing random measure gives (2) and
completes the proof of the Lemma.\bull\enddemo

Let us now define $\mu(s,\tau)\in\Cal M$ by setting $$
\mu(s,\tau)=
\sum_{n=1}^{\infty}\sum_{j=1}^{\infty}a_j(s)a_n(\tau\sigma_js)\delta(
\sigma_n\tau \sigma_is)$$ whenever the series in (1)
converges absolutely and setting $\mu(s,\tau)=0$ otherwise.
It is not difficult to see that the map
$(s,\tau)\to\mu(s,\tau)$ is weak$^*$-Borel.  For a.e.
$\tau$ we have $\mu(s,\tau)=\nu^{S(\tau)}_s$ a.e. and so $$
\int_0^1\|\mu(s,\tau)\|_p^pds \le K_p^p.$$ It follows that
$$ \int_{\Cal T}\int_0^1 \|\mu(s,\tau)\|_p^p ds\,d\tau \le
K_p^p.$$

\proclaim{Lemma 6.6} For almost every $(s,\tau) \in
E\times\Cal T$ we have $$\|\mu(s,\tau)\|_p^p \ge
\left(\sum_{j=1}^M-\sum_{j=M+1}^{\infty}\right)\sum_{n=1}^{\infty}
|a_j(s)|^p|a_n(\tau\sigma_js)|^p.$$\endproclaim

\demo{Proof}Assuming $(s,\tau)$ belongs to set where (1)
holds it is clear the conclusion fails for $(s,\tau)$ if and
only if there exist two distinct pairs $(n,j),(m,i)$ where
$m,n\in\bold N$ and $i,j\le M$ so that
$\sigma_n\tau\sigma_js=\sigma_m\tau\sigma_is$ and
$a_n(\tau\sigma_js)\neq 0.$

Assume then the conclusion of the lemma is false.  Then
there is a distinct pair $(n,j),(m,i)$ as above and Borel
subset $B$ of $E\times \Cal T$ so that $\int_B
|a_n(\tau\sigma_js)|ds\,d\tau >0$ and
$\sigma_n\tau\sigma_js=\sigma_m\tau\sigma_is$ for
$(s,\tau)\in B.$ Note first that we must have $i\neq j.$

It will now follow from Fubini's theorem that there is a
Borel subset $B'$ of $\Gamma$ and a fixed $s\in E$ and
$\pi\in\Pi_P$ for which $\int_{B'}|a_n(\gamma\pi\sigma_i
s)|d\gamma>0 $ and so that
$\sigma_n\gamma\pi\sigma_js=\sigma_m\gamma\pi\sigma_is$ for
$\gamma\in B'.$

Now $\pi\sigma_js\in D(N,k)$ for some $k$ and
$\pi\sigma_is\in D(N,l)$ where $l\neq k$ since $i,j\le M$
and $s\in E$.  We write $\Gamma=\Gamma_k\times \Gamma'$
where $\Gamma'=\prod_{r\neq k}\Gamma_r.$ Again by Fubini's
theorem there exists a fixed $\gamma' \in \Gamma'$ and a
Borel subset $B_0$ of $\Gamma_k$ so that $\int_{B_0}
|a_n(\gamma_k\gamma'\pi\sigma_is)|d\gamma_k>0$ and
$\sigma_n\gamma_k\gamma'\pi\sigma_js=\sigma_m\gamma_k\gamma'\pi\sigma_is$
for $\gamma_k\in B_0.$

Now we note that $\gamma'\pi\sigma_is\in D(N,l)$ is fixed by
every $\gamma_k$ and so
$\sigma_n\gamma_k\gamma'\pi\sigma_js=s'$ is fixed for
$\gamma_k\in B_0.$ But $B_0$ has positive measure and
$\gamma'\pi\sigma_i s\in D(N,k)$.  Thus there is a subset
$A$ of $D(N,k)$ so that $\int |a_n(t)|dt >0$ and
$\sigma_n(A)=\{s'\}.$ This means that $|\nu^T_t|(\{s'\})>0$
on a set of positive measure and we have a
contradiction.\bull\enddemo

We now complete the proof.  We have $$ \align K_p^p &\ge
\int_E\int_{\Cal T}\|\mu(s,\tau)\|_p^p d\tau\, ds \\ &\ge
\int_E\int_{\Cal
T}\left(\sum_{j=1}^{\infty}-2\sum_{j=M+1}^{\infty}\right)
\sum_{n=1}^{\infty}|a_j(s)|^p|a_n(\tau\sigma_js)|^p d\tau\,
ds \endalign $$ As before $$ \int_{\Cal
T}|a_n(\tau\sigma_js)|^pd\tau =\int_0^1|a_n(t)|^pdt.$$ Thus
we obtain $$ K_p^p \ge
\left(\sum_{n=1}^{\infty}|a_n|^pdt\right)\left(\int_E
(\sum_{n=1}^{\infty}-2\sum_{n>M})|a_n(t)|^p dt\right) $$ and
this implies that $$ K_p^p \ge \left(\int_0^1
H_pdt\right)\left(\int_EH_pdt-2\epsilon \right).$$ We
finally deduce that $J=\int_0^1H_pdt$ then
$J(J-3\epsilon)\le K_p^p.$ But $\epsilon>0$ is arbitrary and
so we have $J^2\le K_p^p.$ But as this applies to all such
$T$ we have the conclusion $K_p^2\le K_p$ i.e.  $K_p\le 1.$
Now this applies to all $0<p\le 1.$

Returning to our original $T$ we note that $$ \int_0^1
\sum_{n=1}^{\infty}|a_n(s)|^pds \le 1$$ for all $p$.  Let
$R(s)$ be the number of points in the support of $\nu^T_s.$
Then $$R(s)=\lim_{p\to 0}\sum |a_n(s)|^p.$$ By Fatou's Lemma
$$ \int_0^1 R(s)ds \le 1.$$

To deduce that $T$ is elementary we must show $R(s)=1$ a.e.
If $X$ is order-continuous this is obvious since the fact
that $T$ is surjective requires $R(s)\ge 1$ a.e.

For the general case we again use Proposition 2.1.  We first
note that $R(s)<\infty$ a.e.  We then show that if $x\in X$
then $$ Tx(s)=\sum_{j=1}^{\infty}a_j(s)x(\sigma_js)\tag 3$$
almost everywhere.  To see this it suffices to consider the
case $x\ge 0.$ We first find an increasing sequence of Borel
sets $F_n$ such that $x\chi_{F_n}\in L_{\infty}$ and $\cup
F_n=[0,1].$ Then by Proposition 2.1 $Tx\chi_{F_n}$ converges
in measure to $Tx.$ However, $$Tx\chi_{F_n}(s)
=\sum_{j=1}^{\infty}a_j(s)x(\sigma_js)\chi_{F_n}(\sigma_js)$$
and this converges almost everywhere to the right-hand side
of (3).  Now as in the order-continuous case we can argue
that if $T$ is onto we must have $R(s)\ge 1$ a.e. and hence
$R(s)=1$ a.e.  We conclude that $T$ is elementary in the
case when $X$ has property $(P).$

If $X$ fails property (P) then by Lemma 5.2 $X'$ has
property (P).  Further, Proposition 2.5 says that the
adjoint $T':X'\to X'$ is a surjective isometry.  But then
$T'$ is elementary and by Lemma 2.2 $T''$ and hence $T$ is
elementary.\bull\enddemo

\vskip2truecm

\subheading{7.  Isometries in spaces not isomorphic to
$L_p$}

We now recall the definition of the Boyd indices of an r.i.
space $X$ (cf.  [19] p. 129).  For $0<s<\infty$ define
$D_s:X\to X$ by $D_sf(t) =f(t/s)$ where we let $f(t)=0$ for
$t>1.$ Then the Boyd indices $p_X$ and $q_X$ are defined by
$$ \align \frac1{p_X} &= \lim_{s\to\infty}\frac{\log
\|D_s\|}{\log s}\\ \frac1{q_X}&= \lim_{s\to 0}\frac{\log
\|D_s\|}{\log s}.  \endalign $$

\proclaim{Proposition 7.1}Let $X$ be an r.i. space and
suppose $T:X\to X$ is an elementary operator.  Suppose
$p_X\le r\le q_X.$ Then $T$ is bounded on $L_r[0,1]$ and
$\|T\|_{L_r}\le \|T\|_X.$ \endproclaim

This Proposition is proved in [16] Theorem 5.1.  In fact the
hypotheses of [16] Theorem 5.1 suppose $X$ is a quasi-Banach
space and have an additional unnecessary restriction $r\le
\min(1,q_X).$ This restriction is not used in Theorem 5.1 of
[16] but is important in the following Theorem 5.2.

We will however show a direct proof under the assumption
that $Tx=ax\circ\sigma$ where $\sigma$ is a Borel
automorphism of $[0,1]$ which is the case we need.  For
convenience we consider the case $r<\infty$, the other case
being similar.  Let us assume $\|T\|_X=1.$ We define a
measure Borel measure $\mu$ by
$\mu(B)=\lambda(\sigma^{-1}B)$ and it follows from the fact
that $T$ is bounded that $\mu$ is continuous with respect to
$\lambda$ and so has a Radon-Nikodym derivative $w.$ Now for
any $x$ we have $$ \align \|Tx\|_r^r &= \int_0^1|a(s)|^r
|x(\sigma(s))|^rds\\ &= \int_0^1 |a(\sigma^{-1}s)|^rw(s)
|x(s)|^rds \endalign $$ and so we need to show that
$|a(\sigma^{-1}s)|^rw(s)\le 1$ a.e.  Suppose not.  Then
there is a Borel set $E$ of positive measure $\delta$ and
$0<\alpha,\beta$ so that $\alpha^r\beta>1$ and
$|a(\sigma^{-1}s)|>\alpha$ and $w(s)>\beta$ for $s\in E.$
Then if $x$ is supported in $E$ it quickly follows that
$\|Tx\|_X \ge \alpha \|D_{\beta}x\|_X$ and so
$\|D_{\beta}\|_{X[0,\delta]} \le \alpha^{-1}.$ However for
any $\delta>0$ we have the estimate
$\|D_{\beta}\|_{X[0,\delta]} \ge
\max(\beta^{1/p},\beta^{1/q})$ where $p=p_X$ and $q=q_X.$
Thus $\beta^{1/r} \le \alpha^{-1}$ contrary to
assumption.\bull\enddemo

In [18] Lamperti shows that if $1<p<\infty$ then $L_p[0,1]$
has an equivalent r.i. norm (not equal to the original norm)
so that there are isometries of the form $Tf=af\circ\sigma$
with $|a|\neq 1$ on a set of positive measure.  In the next
theorem we show that if $X$ is not equal to $L_p$ up to
equivalence of norm then the isometries of $X$ can only be
of the very simplest form.

\proclaim{Theorem 7.2}Suppose $X$ is an r.i. space and that
$T$ is a surjective isometry.  Then either $X=L_p[0,1]$ up
to equivalence of norm for some $1\le p\le\infty$ or there
is an invertible measure-preserving Borel map
$\sigma:[0,1]\to [0,1]$ and a function $a\in L_0[0,1]$ with
$|a|=1$ a.e. such that $Tx=ax\circ\sigma$ for $x\in X.$

\demo{Remark} If $p_X<q_X$ this follows routinely from
Proposition 7.1.  The interesting case is thus when
$p_X=q_X.$\enddemo

\demo{Proof}We have that $Tx=ax\circ\sigma$ where $|a|>0$
a.e. and $\sigma:([0,1],\lambda)\to ([0,1],\lambda)$ is a
Borel automorphism by Theorem 6.4.  Suppose $p_X\le r\le
q_X.$ By Proposition 6.1, $\|Tx\|_r\le \|x\|_r$ whenever
$x\in L_r$ and similarly $\|T^{-1}x\|_r\le \|x\|_r.$ Thus
$T$ also defines an isometry on $L_r$ for $p_X\le r\le q_X.$

Let us consider first the case $p_X=q_X=\infty.$ Then
$|a|\le 1$ a.e.  In fact if $B=\{|a|<1-\epsilon\}$ for some
$\epsilon>0$ then $T\chi_{\sigma B}=a\chi_B$ since $\sigma$
is invertible and so $\lambda(B)=0.$ Hence $|a|=1$ a.e.
Again suppose $\epsilon>0$.  Then, assuming $X$ is not
isomorphic to $L_{\infty},$ there exists a least $\delta$ so
that $\|\chi_{[0,\delta]}\|_X=\|\chi_{[0,\epsilon]}\|_X.$
Then if $\lambda(B)=\delta$ we have
$\lambda(\sigma^{-1}(B))\ge \delta.$ For an arbitrary Borel
subset $E$ of $[0,1],$ we can split $E$ into sets of measure
$\delta$ and one remainder set to conclude
$\lambda(\sigma^{-1}(E))\ge \lambda(E)-\delta.$ As
$\epsilon$ was arbitrary $\lambda(\sigma^{-1}(E))\ge
\lambda(E)$ for all $E$.  Since $\sigma$ is invertible this
forces $\lambda(\sigma^{-1}(E))=\lambda(E)$ for every $E$
i.e.  $\sigma$ is measure-preserving.

We turn to the case when $p_X=p<\infty.$ It then follows
that if $|a|=1$ a.e. we must have $\sigma$
measure-preserving.  We thus assume that $|a|\neq 1$ on a
set of positive measure; we will prove that the norm
$\|\,\|_X$ is equivalent to $\|\,\|_p.$ It suffices to
consider the case when $a>0.$ It follows first that
$\{a>1\}$ and $\{a<1\}$ both must have positive measure.

Let us now make an assumption.

\proclaim{Assumption}There exist two disjoint closed
intervals $I_1$ and $I_2$ contained in $(1,\infty)$ and so
that $a^{-1}(I_1)$ and $a^{-1}(I_2)$ have positive
measure.\endproclaim

We will proceed under this assumption.  We can then deduce
that there is a constant $\kappa>1$ and two disjoint Borel
sets $A_1,A_2$ of positive measure such that $a(s)>\kappa$
for $s\in A_2$, while $a(s)>\kappa a(t)$ whenever $s\in
A_1,t\in A_2$ but $a(s)\le \kappa a(t)$ whenever $s,t$ are
either both in $A_1$ or both in $A_2.$

Let
$\delta=\min(\lambda(\sigma(A_1)),\lambda(\sigma(A_2))).$
Let us first note that if $x$ is supported in
$\sigma(A_1)\cup\sigma(A_2)$ then since $a>\kappa$ on
$A_1\cup A_2$ we will have $\lambda(\text{supp }Tx)\le
\kappa^{-p}\lambda(\text{supp }x).$ We also have since $a$
is bounded on $A_1\cup A_2$ an estimate $\lambda(\text{supp
}Tx)\ge c\lambda(\text{supp }x)$ for some $c>0.$

Let us consider any nonnegative $x\in X$ with support $E$ of
measure at most $\delta,$ and such that $\|x\|_p=1.$ We
define the distortion $H(x)$ by setting $$H(x)=\text{ess
sup}\{x(s)/x(t):\ (s,t)\in E^2 \}.$$ If the distortion
$H(x)<\infty$ then it is clear we can define
$\alpha(x)=\text{ess inf}\{x(s):\ s\in E\}$ and $\beta(x)
=\text{ess sup}\{x(s):\ s\in E\}$ and then
$0<\alpha(x)<\beta(x)<\infty$ and $\beta(x)=H(x)\alpha(x).$
Further $\alpha(x)\chi_E\le x\le \beta(x)\chi_E.$

We now define a procedure.  Assume $H(x)<\infty.$ Given such
$x$ we define $x'$ with the same distribution supported on
$\sigma(A_1)\cup\sigma(A_2)$ so that $x'\le
(\alpha(x)\beta(x))^{1/2}$ on $\sigma(A_1)$ but $x'\ge
(\alpha(x)\beta(x))^{1/2}$ on $\sigma(A_2)\cap\text{supp
}x'.$

Now compute $y=Tx'.$ Then $y$ is supported on $A_1\cup A_2.$
If $y(s),y(t)$ are both nonzero and $s,t$ are in the same
$A_j$ we have $y(s)\le \beta(x)^{1/2}\alpha(x)^{-1/2}\kappa
y(t).$ If $s\in A_1$ and $t\in A_2$ we have $y(s) \le
\beta(x)\alpha(x)^{-1}\kappa^{-1}y(t)$.  If $s\in A_2$ and
$t\in A_1$ we have $y(s) \le \kappa y(t).$ It follows that
$$ H(y) \le \max(\kappa H(x)^{1/2},\kappa^{-1}H(x)).$$
Notice also that $c\lambda(\text{supp }x) \le
\lambda(\text{supp }y)\le \kappa^{-p}\lambda(\text{supp
}x).$

If we put $y=x_1$ we can then iterate the procedure to
produce a sequence $(x_n).$ Let $\delta_n=\lambda(\text{supp
}x_n);$ then $c\delta_{n}\le \delta_{n+1}\le
\kappa^{-p}\delta_n$ and, in particular,
$\lim_{n\to\infty}\delta_n=0.$ Since $$ H(x_n) \le
\max(\kappa H(x_{n-1})^{1/2},\kappa^{-1}H(x_{n-1})),$$ we
deduce that $\limsup H(x_n) <\kappa^5.$

Fix any $n$ where $H(x_n)<\kappa^5.$ Then for suitable
$\alpha>0$ and a Borel set $E$ of measure $\delta_n$ we have
$\alpha \chi_E\le x_n\le \kappa^5\alpha\chi_E.$ However
$\|x_n\|_p=1$ and so we obtain $ \alpha\delta_n^{1/p} \le 1
\le \kappa^5\alpha\delta_n^{1/p}$ or $$
\kappa^{-5}\delta_n^{-1/p}\le \alpha \le \delta_n^{-1/p}.$$

Now we introduce the notation $\phi(t)=\|\chi_{[0,t]}\|_X$.
The above inequalities give us $$ \alpha \phi(\delta_n) \le
\|x_n\|_X=\|x\|_X\le \kappa^5\alpha \phi(\delta_n),$$ and
hence $$ \kappa^{-5}\phi(\delta_n)\delta_n^{-1/p} \le\|x\|_X
\le \kappa^5\phi(\delta_n)\delta_n^{-1/p}.$$

Now since $\delta_{n+1}\ge c\delta_n$ we have that for
$\delta_{n+1}\le t\le \delta_n,$

$$ c^{1/p}\kappa^{-5}\phi(t)t^{-1/p}\le \|x\|_X \le
c^{-1/p}\kappa^5\phi(t)t^{-1/p}.$$

As $H(x_n)\le \kappa^5$ eventually we can conclude that
$$0<\liminf_{t\to 0}\phi(t)t^{-1/p} <\limsup_{t\to 0}
\phi(t)t^{-1/p} <\infty.$$ In fact if we let $K=\limsup
\phi(t)t^{-1/p}$ we obtain $$ c\kappa^{-5}K\le \|x\|_X \le
c^{-1}\kappa^5K.$$

But this estimate is independent of the original choice of
$x$ subject to $\lambda(\text{supp }x)\le \delta$,
$H(x)<\infty$ and $\|x\|_p=1.$ Hence we obtain that
$\|x\|_X$ is equivalent to $\|x\|_p.$

Thus our assumption yields the conclusion that $X=L_p[0,1]$
up to an equivalent norm.  Clearly it suffices to find one
surjective isometry for which the assumption holds to give
this conclusion.

If the assumption fails for $T$ then $a$ is essentially
constant (with value $\alpha$, say) on $\{a>1\}.$ If the
assumption fails for $T^{-1}$ it is easy to see that $a$ is
also essentially constant (with value $\beta,$ say) on
$\{a<1\}.$ Now the same reasoning must apply to any
surjective isometry.  However it is now easy to construct an
isometry of the form $S=TV_{\tau_1}TV_{\tau_2}T$, where
$V_{\tau}=x\circ\tau$ for some measure preserving Borel
automorphism $\tau,$ and so that $S\chi_{[0,1]}$ takes each
of the four distinct values
$\alpha^3,\alpha^2\beta,\alpha\beta^2$ and $\beta^3$ (of
which three must be distinct from 1) with positive measure.
Thus we can again conclude that $X$ is isomorphic to
$L_p.\bull$\enddemo

\demo{Remarks} This theorem can be cast as a statement about
maximal norms (cf.  [25], [17]).  A Banach space $X$ has a
{\it maximal norm} if no equivalent norm has a strictly
bigger group of invertible isometries.  The above theorem
shows immediately that any r.i. space on $[0,1]$ which is
not isomorphic to $L_p[0,1]$ has a maximal norm; Rolewicz
[25] showed that the spaces $L_p[0,1]$ have maximal norms.
However if $X$ is isomorphic but not isometric to $L_p$ its
norm cannot be maximal; this follows rather easily from
Proposition 7.1 and the almost transitivity of the norm in
$L_p.$ \enddemo

\vskip2truecm

\subheading{References}

\item{1.} S. Banach, {\it Theorie des operations lineaires,}
Warsaw 1932.

\item{2.} F.F.  Bonsall and J. Duncan, {\it Numerical ranges
of operators on normed spaces and of elements of normed
algebras,} London Math.  Soc.  Lecture Notes, 2, Cambridge
Univ.  Press, 1971.

\item{3.} N.L.  Carothers, S.J.  Dilworth and D.A.
Trautman, On the geometry of the unit sphere of the Lorentz
space $L_{w,1}$, Glasgow Math.  J. 34 (1992) 21-27.

\item{4.} N.L.  Carothers, R.G.  Haydon and P.K.  Lin, On
the isometries of the Lorentz space $L_{w,p},$ preprint

\item{5.} N.L.  Carothers and B. Turett, Isometries of
$L_{p,1},$ Trans.  Amer.  Math.  Soc. 297 (1986) 95-103.

\item{6.} R.J.  Fleming, J. Goldstein and J.E.  Jamison, One
parameter groups of isometries on certain Banach spaces,
Pacific J. Math. 64 (1976) 145-151.

\item{7.} R.J.  Fleming and J.E.  Jamison, Isometries of
certain Banach spaces, J. London Math.  Soc. 9 (1974)
121-127.

\item{8.} R.J.  Fleming and J.E.  Jamison, Isometries of
Banach spaces - a survey, preprint.

\item{9.} G. Godefroy and N.J.  Kalton, The ball topology
and its applications, Contemporary Math. 85 (1989) 195-238.

\item{10.} J. Goldstein, Groups of isometries on Orlicz
spaces, Pacific J. Math. 48 (1973) 387-393.

\item{11.} R. Grzaslewicz, Isometries of $L^1\cap L^p,$
Proc.  Amer.  Math.  Soc. 93 (1985) 493-496.

\item{12.} R. Grzaslewicz and H.H.  Schaefer, Surjective
isometries of $L^1\cap L^{\infty}[0,\infty)$ and
$L^1+L^{\infty}[0,\infty),$ Indag.  Math. 3 (1992) 173-178.

\item{13.} J. Jamison, A. Kaminska and P.K.  Lin, Isometries
in some classes of generalized Orlicz spaces, preprint.

\item{14.} N.J.  Kalton, The endomorphisms of $L_p,\ 0\le
p\le 1,$ Indiana Univ.  Math.  J. 27 (1978) 353-381.

\item{15.} N.J.  Kalton, Representations of operators on
function spaces, Indiana Univ.  Math.  J. 33 (1984) 640-665.

\item{16.} N.J.  Kalton, Endomorphisms of symmetric function
spaces, Indiana Univ.  Math.  J. 34 (1985) 225-247.

\item{17.} N.J.  Kalton and G.V.  Wood, Orthonormal systems
in Banach spaces and their applications, Math.  Proc.  Camb.
Phil.  Soc. 79 (1976) 493-510.

\item{18.} J. Lamperti, On the isometries of some function
spaces, Pacific J. Math. 8 (1958) 459-466.

\item{19.} J. Lindenstrauss and L. Tzafriri, {\it Classical
Banach spaces, Vol. 2, Function spaces,} Springer 1979.

\item{20.} G. Lumer, Semi-inner product spaces, Trans.
Amer.  Math.  Soc. 100 (1961) 29-43.

\item{21.} G. Lumer, Isometries of reflexive Orlicz spaces,
Bull.  Amer.  Math.  Soc. 68 (1962) 28-30.

\item{22.} G. Lumer, On the isometries of reflexive Orlicz
spaces, Ann.  Inst.  Fourier 13 (1963) 99-109.

\item{23.} G. Lumer and R.S.  Phillips, Dissipative
operators in a Banach space, Pacific J. Math. 11 (1961)
679-698.

\item{24.} D.M.  Oberlin, Translation-invariant operators on
$L^p(G),\ 0<p<1,$ Michigan Math.  J. 23 (1976) 119-122.

\item{25.} S. Rolewicz, {\it Metric linear spaces,} Polish
Scientific Publishers, Warsaw, 1972.

\item{26.} H.P.  Rosenthal, Contractively complemented
subspaces of Banach spaces with reverse monotone
(transfinite) bases, Longhorn Notes, The University of Texas
Functional Analysis seminar, 1984-5, 1-14.

\item{27.} W. Rudin, {\it Functional Analysis, Second
Edition,} McGraw-Hill, New York, 1991.

\item{28.} A.R.  Sourour, Pseudo-integral operators, Trans.
Amer.  Math.  Soc. 253 (1979) 339-363.

\item{29.} L. Weis, On the representation of
order-continuous operators by random measures, Trans.  Amer.
Math.  Soc. 285 (1984) 535-564.

\item{30.} P. Wojtaszczyk, {\it Banach spaces for analysts,}
Cambridge University Press 1991.

\item{31.} M.G. Zaidenberg, Groups of isometries of Orlicz
spaces, Soviet Math.  Dokl. 17 (1976) 432-436.

\item{32.} M.G. Zaidenberg, On the isometric classification
of symmetric spaces, Soviet Math.  Dokl. 18 (1977) 636-640.

\item{33.} M.G. Zaidenberg, Special representations of
isometries of functional spaces, Investigations on the
theory of functions of several real variables, Yaroslavl,
1980 (Russian).

\bye